\documentclass[12pt]{article}
\usepackage{amsmath}
\usepackage{amsfonts}
\usepackage{amssymb}
\usepackage{theorem}

\title{Locally compact groups and continuous logic}
\author{A.Ivanov 
\thanks{The research is supported by Polish National Science Centre grant DEC2011/01/B/ST1/01406} 
}

\newtheorem{theorem}{Theorem}[section]
\newtheorem{proposition}[theorem]{Proposition}
\newtheorem{corollary}[theorem]{Corollary}
\newtheorem{lemma}[theorem]{Lemma}
\newtheorem{definition}[theorem]{Definition}

\newtheorem{remark}[theorem]{Remark}

\begin{document} 
\topmargin = 12pt
\textheight = 630pt 
\footskip = 39pt 

\maketitle

\begin{quote}
{\bf Abstract}
We study expressive power of continuous logic in classes of (locally compact) groups. 
We also describe locally compact groups which are separably categorical structures.
\end{quote}

\section{Introduction}

In this paper we give very concrete applications of 
continuous logic in group theory. 
We consider classes of (locally compact) metric groups 
which can be also viewed as (reducts of) axiomatizable 
classes of continuous structures. 
Then by some standard logical tricks we obtain 
several interesting consequences. 
Usually we concentrate on classes which are 
typical in geometric group theory.  \parskip0pt 

The following notion is one of the main objects of the paper. 
A class of groups $\mathcal{K}$ is called {\em bountiful} if 
for any pair of infinite groups $G\le H$ with $H\in \mathcal{K}$ 
there is $K\in \mathcal{K}$ such that $G\le K\le H$ and $|G|=|K|$. 
It was introduced by Ph.Hall and was studied in 
papers \cite{KM}, \cite{phillips}, \cite{sabbagh} and \cite{thomas}. 
Some easy logical observations from \cite{KM} show that if 
$\mathcal{K}$ is a reduct of a class axiomatizable in 
$L_{\omega_1 \omega}$ then $\mathcal{K}$ is bountiful. 
\parskip0pt 

When one considers topological groups, the definition of 
bountiful classes should be modified as follows. 
\begin{definition} 
A class of topological groups $\mathcal{K}$ is called  bountiful 
if for any pair of infinite groups $G\le H$ with $H\in \mathcal{K}$ 
there is $K\in \mathcal{K}$ such that $G\le K\le H$ and 
the density character of $G$ 
(i.e. the smallest cardinality of a dense subset of the space) 
coincides with the density character of $K$. 
\end{definition} 
We mention paper \cite{HHM} where similar questions were studied 
in the case of locally compact groups. 
We will see below that under additional assumptions of metricity 
logical tools become helpful in this class of groups.  
We should only replace first-order logic (or $L_{\omega_1 \omega}$) 
by continuous one. 
We concentrate on negations of properties {\bf (T)}, {\bf FH}, 
{\bf F}$\mathbb{R}$ (\cite{BHV}, \cite{HV}) and on negations 
of boundedness properties classified in \cite{rosendalN}. 

In the final part of the paper we consider separable locally compact groups 
which have separably categorical continuous theory, i.e. the group 
is determined uniquely (up to metric isomorphism) by its continuous theory 
and the the density character. 
It is interesting that some basic properties of the automorphism groups 
of such structures are strongly connected with some classes examined on 
bountifulness below.

In the rest of this introduction we briefly remind 
the reader some preliminaries of continuous logic. 
Then we finish this section by some remarks on sofic groups. 

\bigskip

\paragraph{Continuous structures.} 

We fix a countable continuous signature 
$$
L=\{ d,R_1 ,...,R_k ,..., F_1 ,..., F_l ,...\}. 
$$  
Let us recall that a {\em metric $L$-structure} 
is a complete metric space $(M,d)$ with $d$ bounded by 1, 
along with a family of uniformly continuous operations on $M$ 
and a family of predicates $R_i$, i.e. uniformly continuous maps 
from appropriate $M^{k_i}$ to $[0,1]$.   
It is usually assumed that to a predicate symbol $R_i$ 
a continuity modulus $\gamma_i$ is assigned so that when 
$d(x_j ,x'_j ) <\gamma_i (\varepsilon )$ with $1\le j\le k_i$ 
the corresponding predicate of $M$ satisfies 
$$ 
|R_i (x_1 ,...,x_j ,...,x_{k_i}) - R_i (x_1 ,...,x'_j ,...,x_{k_i})| < \varepsilon . 
$$ 
It happens very often that $\gamma_i$ coincides with $id$. 
In this case we do not mention the appropriate modulus. 
We also fix continuity moduli for functional symbols. 
Note that each countable structure can be considered 
as a complete metric structure with the discrete $\{ 0,1\}$-metric. 

By completeness continuous substructures of a continuous structure are always closed subsets. 

Atomic formulas are the expressions of the form $R_i (t_1 ,...,t_r )$, 
$d(t_1 ,t_2 )$, where $t_i$ are terms (built from functional $L$-symbols). 
In metric structures they can take any value from $[0,1]$.   
{\em Statements} concerning metric structures are usually 
formulated in the form 
$$
\phi = 0 
$$ 
(called an $L$-{\em condition}), where $\phi$ is a {\em formula}, 
i.e. an expression built from 
0,1 and atomic formulas by applications of the following functions: 
$$ 
x/2  \mbox{ , } x\dot- y= max (x-y,0) \mbox{ , } min(x ,y )  \mbox{ , } max(x ,y )
\mbox{ , } |x-y| \mbox{ , } 
$$ 
$$ 
\neg (x) =1-x \mbox{ , } x\dot+ y= min(x+y, 1) \mbox{ , } sup_x \mbox{ and } inf_x . 
$$ 
A {\em theory} is a set of $L$-conditions without free variables 
(here $sup_x$ and $inf_x$ play the role of quantifiers). 
   
It is worth noting that any formula is a $\gamma$-uniformly continuous 
function from the appropriate power of $M$ to $[0,1]$, 
where $\gamma$ is the minimum of continuity moduli of $L$-symbols 
appearing in the formula. 

The condition that the metric is bounded by $1$ is not necessary. 
It is often assumed that $d$ is bounded by some rational number $d_0$. 
In this case the (dotted) functions above are appropriately modified.  
Sometimes predicates of continuous structures map $M^n$ to some 
$[q_1 ,q_2 ]$ where $q_1 ,q_2 \in \mathbb{Q}$.  

The following theorem is one of the main tools of this paper. 
\bigskip 
 
{\bf L\"{o}wenheim-Skolem Theorem.} (\cite{BYBHU}, Proposition 7.3) 
{\em Let $\kappa$ be an infinite cardinal number and assume 
$|L|\le \kappa$. 
Let $M$ be an $L$-structure and suppose $A\subset M$ has 
density $\le\kappa$. 
Then there exists a substructure $N\subseteq M$ containing $A$ such that 
$density(N) \le\kappa$ and $N$ is an elementary substructure of $M$, i.e. 
for every $L$-formula $\phi (x_1 ,...,x_n)$ and $a_1 ,...,a_n \in N$ 
the values of $\phi (a_1 ,...,a_n )$ in $N$ and in $M$ are the same.}  

\bigskip 

\begin{remark} \label{Hofmann}
{\em 
It is proved in \cite{HHM} that for any locally compact group $G$, 
the entire interval of cardinalities between $\aleph_0$ and $w(G)$, 
the weight of the group, is occupied by the weights of closed 
subgroups of $G$.   
We remind the reader that the weight of a topological space $(X,\tau )$ 
is the smallest cardinality which can be realized as 
the cardinality of a basis of $(X,\tau )$.  
If the group $G$ is metric, the weight of $G$ coincides 
with the density character of $G$. 
This yelds the following version of the L\"{o}wenheim-Skolem Theorem. }

 Let $G$ be a locally compact group which is a continuous structure. 
Then for any cardinality $\kappa < density(G)$ there is 
a closed subgroup $H<G$ such that $density(H)=\kappa$ 
and $H$ is an elementary substructure of $G$. 
\end{remark} 

\begin{remark} \label{FaHaSh} 
{\em 
Following Section 4.2 of \cite{fhs} we define a topology 
on $L$-formulas relative to a given continuous theory $T$. 
For $n$-ary formulas $\phi$ and $\psi$ of the same sort set 
$$ 
{\bf d}^T_{\bar{x}} (\phi ,\psi) = sup \{ |\phi (\bar{a}) -\psi (\bar{a} )|: \bar{a} \in M, M\models T\} . 
$$ 
The function ${\bf d}^T_{\bar{x}}$ is a pseudometric. 
The language $L$ is called {\em separable} if 
for every $L$-theory $T$ and any tuple $\bar{x}$ 
the density character of ${\bf d}^T_{\bar{x}}$ is countable. 
By Proposition 4.5 of \cite{fhs} in this case for every 
$L$-model $M$ the set of all interpretations of $L$-formulas 
in $M$ is separable in the uniform topology. 
By Corollary 4.7 of \cite{fhs} if in the formulation of 
the L\"{o}wenheim-Skolem theorem we replace the assumption 
$|L|\le \kappa$ by the condition that $L$ is separable 
then the statement also holds for $\kappa =\aleph_0$. 
This will be applied in Section 2.  
}
\end{remark}

Definability in continuous structures is introduced as follows. 

\begin{definition} 
Let $A\subseteq M$. 
A predicate $P:M^n \rightarrow [0,1]$ is definable in $M$ over $A$
if there is a sequence $(\phi_k (x) :k\ge 1 )$ of $L(A)$-formulas 
such that predicates interpreting $\phi_k (x)$ in $M$ 
converge to $P(x)$ uniformly in $M^n$. 
\end{definition} 

We define the automorphism group $Aut(M)$ of $M$ to be the subgroup 
of $Iso (M,d)$ consisting of all isometries preserving 
the values of atomic formulas. 
It is easy to see that $Aut(M)$ is a closed subgroup with 
respect to the pointwise convergence topology on $Iso (M,d)$. 

The following statement is Corollary 9.11 of \cite{BYBHU}. 
\begin{quote} 
Let $M$ be an $L$-structure with $A\subseteq M$ and suppose 
$P:M^n \rightarrow [0,1]$ is a predicate. 
Then $P$ is definable in $M$ over $A$ if and only if 
whenever $(N,Q)$ is an elementary extension of $(M,P)$, 
the predicate $Q$ is invariant under all automorphisms of $N$ 
that leave $A$ fixed pointwise. 
\end{quote}  

A tuple $\bar{a}$ from $M^n$ is {\em algebraic} in $M$ over 
$A$ if there is a compact subset $C\subseteq M^n$ such that 
$\bar{a}\in C$ and the distance predicate $dist(\bar{x},C)$ 
is definable in $M$ over $A$. 
Let $acl(A)$ be the set of all elements algebraic over $A$. 
In continuous logic the concept of algebraicity is 
parallel to that in traditional model theory 
(see Section 10 of \cite{BYBHU}).  

For every $c_1 ,...,c_n \in M$ and $A\subseteq M$ 
we define the $n$-type $tp(\bar{c}/A)$ of $\bar{c}$ over $A$ 
as the set of all $\bar{x}$-conditions with parameters from $A$ 
which are satisfied by $\bar{c}$ in $M$.  
Let $S_n (T_A )$ be the set of all $n$-types over $A$ 
of the expansion of the theory $T$ by constants from $A$. 
There are two natural topologies on this set. 
The {\em logic topology} is defined by the basis consisting of 
sets of types of the form $[\phi (\bar{x})<\varepsilon ]$, 
i.e. types containing some $\phi (\bar{x})\le \varepsilon'$ with 
$\varepsilon '<\varepsilon$.   
The logic topology is compact. 

The $d$-topology is defined by the metric 
$$
d(p,q)= inf \{ max_{i\le n} d(c_i ,b_i )| \mbox{ there is a model } M \mbox{ with } M\models p(\bar{c})\wedge q(\bar{b})\}. 
$$ 
By Propositions 8.7 and 8.8 of \cite{BYBHU} the $d$-topology is finer 
than the logic topology and $(S_n (T_A ),d)$ is a complete space.

\paragraph{Separable categoricity.} 

A theory $T$ is {\em separably categorical} if any 
two separable models of $T$ are isomorphic. 
By Theorem 12.10 of \cite{BYBHU} a complete theory $T$ is separably 
categorical if and only if for each $n>0$, every $n$-type $p$ is principal. 
The latter means that for every model $M\models T$, the predicate 
$dist(\bar{x},p(M))$ is definable over $\emptyset$. 

Another property equivalent to separable categoricity states that 
for each $n>0$, the metric space $(S_n (T),d)$ is compact.  
In particular for every $n$ and every $\varepsilon$ there is 
a finite family of principal $n$-types $p_1 ,...,p_m$ so that 
their $\varepsilon$-neighbourhoods cover $S_n(T)$. 

In first order logic a countable structure $M$ is 
$\omega$-categorical if and only if $Aut(M)$ is 
an {\em oligomorphic} permutation group, i.e. 
for every $n$, $Aut(M)$ has finitely many orbits 
on $M^n$. 
In continuous logic we have the following modification.  

\begin{definition} 
An isometric action of a group $G$ on a metric space $({\bf X},d)$ 
is said to be approximately oligomorphic if for every $n\ge 1$ and $\varepsilon >0$ 
there is a finite set $F\subset {\bf X}^n$ such that 
$$ 
G\cdot F = \{ g\bar{x} : g\in G \mbox{ and } \bar{x}\in F\}
$$
is $\varepsilon$-dense in $({\bf X}^n,d)$. 
\end{definition} 

Assuming that $G$ is the automorphism group of a non-compact 
separable continuous metric structure $M$, $G$ is approximately 
oligomorphic if and only if  the structure $M$ is separably categorical 
(C. Ward Henson, see Theorem 4.25 in \cite{scho}). 
It is also known that separably categorical structures are 
{\em approximately homogeneous} in the following sense: 
if $n$-tuples $\bar{a}$ and $\bar{c}$ have the same types 
(i.e. the same values $\phi (\bar{a})=\phi (\bar{b})$ for all $L$-formulas $\phi$) 
then for every $c_{n+1}$ and $\varepsilon >0$ there is 
an tuple $b_1 ,...,b_n ,b_{n+1}$ of the same type with 
$\bar{c},c_{n+1}$, so that $d(a_i, b_i )\le \varepsilon$ for $i\le n$.  
In fact for any $n$-tuples $\bar{a}$ and $\bar{b}$ there is 
an automorphism $\alpha$ of $M$ such that 
$$
d(\alpha (\bar{c}),\bar{a})\le d(tp(\bar{a}),tp(\bar{c})) +\varepsilon . 
$$  
(i.e $M$ is {\em strongly $\omega$-near-homogeneous} in the sense
of Corollary 12.11 of \cite{BYBHU}). 

\begin{definition} 
A topological group $G$ is called Roelcke precompact 
if for every open neighborhood of the identity $U$, 
there exists a finite subset $F\subset G$ such that $G=UFU$. 
\end{definition} 

The following theorem is a combination of the remark above, 
Theorem 6.2 of \cite{rosendal},  
Theorem 2.4 of \cite{tsankov} and Proposition 1.20 of \cite{rosendalN}. 

\begin{theorem} 
Let $G$ be the automorphism group of a non-compact separable structure $M$. 

Then \\ 
(i) the group $G$ is approximately oligomorphic if and only if $M$ is separably categorical; \\  
(ii) if $G$ is Roelcke precompact and approximately oligomorphic for 1-orbits, 
then $M$ is separably categorical;\\  
(iii) if the structure $M$ is separably categorical, then $G$ is Roelcke precompact. 
\end{theorem}

\paragraph{Axiomatizability in continuous logic, topological properties and  sofic groups. } 

Suppose $\mathcal{C}$ is a class of metric $L$-structures. 
Let $Th^c (\mathcal{C})$ be the set of all closed $L$-conditions 
which hold in all structures of $\mathcal{C}$. 
It is proved in \cite{BYBHU} (Proposition 5.14 and Remark 5.15) 
that every model of $Th^c (\mathcal{C})$ is elementary equivalent 
to some ultraproduct of structures from $\mathcal{C}$. 
Moreover by Proposition 5.15 of \cite{BYBHU} we have the following statement. 
\begin{quote} 
The class $\mathcal{C}$ is axiomatizable in continuous logic 
if an only if it is closed under metric isomorphisms and 
ultraproducts and its complement is closed under ultrapowers. 
\end{quote} 
Let $Th^c_{\sup} (\mathcal{C})$ be the set of all closed 
$L$-conditions of the form 
$$ 
sup_{x_1} sup_{x_2} ... sup_{x_n} \varphi =0 
\mbox{ ( $\varphi$ does not contain $inf_{x_i}$, $sup_{x_i}$ ), } 
$$
which hold in all structures of $\mathcal{C}$.
Some standard arguments also give the following theorem. 

\begin{theorem} \label{axiom} 
(1) The class $\mathcal{C}$ is axiomatizable in continuous logic 
if an only if it is closed under metric isomorphisms, ultraproducts and 
taking elementary submodels. \parskip0pt  

(2) The class $\mathcal{C}$ is axiomatizable in continuous logic 
by $Th^c_{sup} (\mathcal{C})$ if an only if it is closed under 
metric isomorphisms, ultraproducts and taking substructures. 
\end{theorem} 

It is worth noting that when one considers classes axiomatizable 
in continuous logic it is obviously assumed that all operations 
and predicates are uniformly continuous. 
This shows that some topological properties cannot be described 
(axiomatized) in continuous logic. 

Some other obstacles arise from the fact that existentional 
quantifiers cannot be expressed in continuous logic. 
For example consider the class of all 
metric groups which are discrete in their metrics (with $id$ as continuity moduli). 
This class is not closed under metric ultraproducts but 
if we replace all metrics by the $\{ 0,1\}$-one 
we just obtain the (axiomatizable) class of all groups. 

It may also happen that when we extend an axiomatizable 
class of structures with the $\{ 0,1\}$-metric 
\footnote{in this case axiomatizability in continuous logic is equivalent to axiomatizability in first-order logic} 
by (abstract) structures from this class with all possible 
(not only possible discrete) metrics we lose axiomatizability. 
A nice example of this situation is the class of non-abelian groups 
with $[0,1]$-metrics. 
For example there is a sequence of non-abelian groups $G_n \le Sym (2^n +3 )$ 
with $G_n \cong \mathbb{Z}(2)^n \times S_3$ so that their metric unltraproduct 
 with respect to Hamming metrics is abelian (an easy exercise). 

Continuous axiomatizability appears in one of the most active areas in group theory 
as follows. 
\begin{quote} 
An abstract group is {\em sofic} if it is embeddable into a metric 
ultraproduct of finite symmetric groups with Hamming metrics. 
\end{quote} 
Let $\mathcal{S}$ be the class of complete $id$-continuous 
metric groups of diameter 1, which are embeddable as closed subgroups  
via isometric morphisms into a metric ultraproduct of finite symmetric 
groups with Hamming metrics.  
This class is axiomatizable by Theorem \ref{axiom}. 
We call it the class of {\em metric sofic groups}. 

\begin{corollary} 
The class of metric sofic groups is $sup$-axiomatizable (i.e. by its theory $Th^c_{sup}$).  
\end{corollary} 

It is folklore that any abstract sofic group can be 
embedded into a metric ultraproduct of finite symmetric 
groups as a discrete subgroup 
(see the proof of Theorem 3.5 of \cite{pestov}). 
This means that the set of all abstract sofic groups consists 
of all discrete structures of the class $\mathcal{S}$.

\section{Boundedness properties} 

It is worth noting that many classes from geometric 
group theory  are just universal. 
For example if a group has free isometric actions 
on real trees (resp. Hilbert spaces) then any its subgroup 
has the same property. 
Similarly a closed subgroup of  a locally compact  amenable 
group is amenable. 
\footnote{the class of discrete initially amenable groups (see \cite{cornulier})  is  universal too} 
Thus these classes are bountiful. 

On the other hand if we extend these classes by non-compact 
locally compact groups without Kazhdan's property ${\bf (T)}$
or by groups admitting isometric actions on real trees  
without fixed points then we lose universality.  
Are these classes still bountiful? 
We may further extend our classes by so called 
non-{\em boundedness properties} introduced in \cite{rosendalN}. 
For example consider metric groups which satisfy non-${\bf OB}$ 
(in terms of \cite{rosendalN}): they have isometric strongly 
continuous actions (i.e. the map $g\rightarrow g\cdot x$ defined on 
$G$ is continuous for each $x$) on metric spaces with unbounded orbits. 
The first part of this section is devoted to some modifications of this property.  
We will show how continuous logic can work in these cases. 
In fact metric groups from these classes can be presented 
as reducts of continuous metric structures  
which induce some special actions.  

In the second part of the section we consider non-${\bf (T)}$ and 
non-${\bf F\mathbb{R}}$ (of fixed points for isometric actions on real trees).    
Note that ${\bf (T)}$ and property ${\bf FH}$ 
(that any strongly continuous isometric affine action on a real Hilbert space has a fixed point)  
are equivalent for $\sigma$-compact locally compact groups 
(see Chapter 2 in \cite{BHV}).  
Since definitions of these properties require Hilbert spaces (or unbounded trees), 
we will here apply a many-sorted version of continuous logic 
(as in Section 15 of \cite{BYBHU}). 
We will also present our groups as a union of an increasing 
chain of subsets of bounded diameters treating each subset as a sort.   
This situation is very natural if the group is $\sigma$-compact 
(i.e. a union of an increasing chain of compact subsets). 

It is worth noting that by Section 1.10 of \cite{rosendalN} 
in the case of $\sigma$-locally compact groups 
(=$\sigma$-compact locally compact) Roelcke precompatness coincides 
with all boundedness properties studied in \cite{rosendalN} excluding only ${\bf FH}$.  
In particular it coincides with compactness and property ${\bf OB}$. 
On the other hand an elementary submodel of a non-compact (resp. compact) continuous 
structure is also non-compact (resp. compact, see \cite{BYBHU}, Section 10). 
Thus by the L\"{o}wenheim-Skolem theorem (in a 1-sorted language) 
non-compacness is bountiful in the class of locally compact groups.  
When the property non-${\bf OB}$ coincides with non-compacness 
(as in the case of locally compact Polish groups) it is also bountiful. 
This explains why in the first part of the section we do not 
assume that a group is locally compact or Polish.

It is worth noting that our methods do not work 
for the classes of (locally compact) groups satisfying 
properties {\bf (T)}, {\bf F}$\mathbb{R}$ and {\bf FH} 
(see discussion before Proposition \ref{discr}). 
The case of amenable Polish groups is open and  looks very intresting. 
A topological group $G$ is called {\em amenable} 
if every $G$-flow admits an invariant Borel probability measure. 
In the case of locally compact groups this definition coincides with 
the classical one. 
It is noticed in \cite{kechrisN}, that the group $Sym(\omega )$ 
of all permutations of $\omega$ is amenable. 
Since it has closed non-amenable subgroups, the class of amenable 
Polish groups is not universal (with respect to taking closed subgroups).

\subsection{Negations of strong boundedness and OB} 

An abstract group $G$ is {\em Cayley bounded} if for every generating subset
$U\subset G$ there exists $n\in \omega$ such that every element
of $G$ is a product of $n$ elements of $U\cup U^{-1}\cup\{ 1\}$.
If $G$ is a Polish group then $G$ is {\em topologically Cayley bounded} 
if for every analytic generating subset $U\subset G$ 
there exists $n\in \omega$ such that every element
of $G$ is a product of $n$ elements of $U\cup U^{-1}\cup\{ 1\}$.
It is proved in \cite{rosendal} that for Polish groups property 
{\bf OB} is equivalent to topological Cayley boundedness together 
with {\em uncountable topological cofinality}: $G$ is not the union 
of a chain of proper open subgroups. 

\paragraph{Discrete groups.} 

Let us consider the abstract (discrete) case. 
A group is {\em strongly bounded} if it is Cayley bounded and
cannot be presented as the union of a strictly increasing chain
$\{ H_n :n\in \omega\}$ of proper subgroups
(has {\em cofinality} $>\omega$).  
It is known that strongly bounded groups have property {\bf FA}, 
i.e. any action on a simplicial tree fixes a point. 

The class of strongly bounded groups is not bountiful. 
Indeed, by \cite{dC} for any finite perfect group $F$ and 
an infinite $I$ the power $F^I$ is strongly bounded. 
Since $F^I$ is locally finite, any its countable subgroup 
has cofinality $\omega$. 
Similar arguments can be applied to property ${\bf FA}$.  

It is shown in \cite{dC}, that strongly bounded groups
have property {\bf FH}.   
It can be also deduced from \cite{dC} that strongly bounded groups 
have property ${\bf F}\mathbb{R}$ that every isometric action of $G$ 
on a real tree has a fixed point (since such a group acting on 
a real tree has a bounded orbit, all the elements are elliptic 
and it remains to apply cofinality $>\omega$). 
It is now clear that the bountiful class of groups having 
free isometric actions on real trees (or on real Hilbert spaces) 
is disjoint from strong boundedness.

\begin{proposition} \label{discr} 
The following classes of groups are reducts of axiomatizable 
classes in $L_{\omega_1 \omega}$: \\ 
(1) The complement of the class of strongly bounded groups; \\ 
(2) The class of groups of cofinality $\le \omega$; \\  
(3) The class of groups which are not Cayley bounded; \\ 
(4) The class of groups presented as non-trivial free products with amalgamation 
(or HNN-extensions); \\   
(5) The class of groups having homomorphisms onto $\mathbb{Z}$.  

All these classes are bountiful. 
The class of groups which do not have property {\bf FA} is bountiful too. 
\end{proposition} 

{\em Proof.} 
(1) We use the following characterization of strongly bounded 
groups from \cite{dC}.  
\begin{quote} 
A group is strongly bounded if and only if for every 
presentation of $G$ as $G=\bigcup_{n\in \omega} X_n$ for 
an increasing sequence $X_n$, $n\in \omega$, with 
$\{ 1\} \cup X^{-1}_n \cup X_n \cdot X_n \subset X_{n+1}$ 
there is a number $n$ such that $X_n =G$. 
\end{quote} 
Let us consider the class $\mathcal{K}_{nb}$ of all structures 
$\langle G, X_n \rangle_{n\in\omega}$ with the axioms 
stating that $G$ is a group, $\{ X_n \}$ is a sequence of 
unary predicates on $G$ defining a strictly increasing 
sequence of subsets of $G$ with 
$\{ 1\} \cup X^{-1}_n \cup X_n \cdot X_n \subset X_{n+1}$ 
(these axioms are first-order) and 
$$ 
(\forall x) (\bigvee_{n\in\omega} x\in X_n ). 
$$
By the L\"{o}wenheim-Skolem theorem for countable fragments of 
$L_{\omega_1 \omega}$ (\cite{keisler}, p.69)
any subset $C$ of such a structure is contained in an elementary 
submodel of cardinality $|C|$ (the countable fragment which we consider 
is the minimal fragment containing our axioms). 
This proves bountifulness in case (1). \\ 
(2) The case groups of cofinality $\le \omega$ is similar. \\ 
(3) The class of groups which are not Cayley bounded is 
a class of reducts of all groups expanded by an unary predicate  
$\langle G,U\rangle$ with an $L_{\omega_1 \omega}$-axiom stating 
that $U$ generates $G$ and with a system of first-order axioms 
stating that there exists an element of $G$ which is not a 
product of $n$ elements of $U\cup U^{-1} \cup \{ 1 \}$. 
The rest is clear.   \\ 
(4) The class of groups which can be presented as 
non-trivial free products with amalgamation is the class 
of reducts of all groups expanded by two unary predicates  
$\langle G,U_1 ,U_2 \rangle$ with first-order axioms that $U_1$ and $U_2$ 
are subgroups and with $L_{\omega_1 \omega}$-axioms 
stating tha $U_1 \cup U_2$ generates $G$ and a word in the 
alphabeth  $U_1 \cup U_2$ is equal to 1 if and only if this 
word follows from the relators of the free product of $U_1$ 
and $U_2$ amalgamated over $U_1 \cap U_2$. 
The rest of (4) is clear.  \\ 
(5) Groups having homomorphisms onto $\mathbb{Z}$ 
can be considered as reducts of structures in the language 
$\langle \cdot ,...U_{-n},...,U_0 ,...,U_{m},...\rangle$, 
where predicates $U_t$ denote preimages of the corresponding integer numbers.

To see that the class of groups without {\bf FA} is 
bountiful, take any infinite $G\models {\bf notFA}$.  
It is well-known (\cite{serre}, Section 6.1) 
that such a group belongs to the union of 
the classes from statements (2),(4) and (5).  
Thus $G$ has an expansion as in one of the cases (2),(4) or (5). 
Now applying the L\"{o}wenheim-Skolem theorem, for any  
$C\subset G$ we find a subgroup of $G$ of cardinality $|C|$ 
which contains $C$ and does not satisfy {\bf FA}.   
$\Box$

\bigskip 


\paragraph{Topological groups.} 

As we already mentioned in Introduction separably categorical structures have 
Roelcke precompact automorphism groups. 
In the following definition we consider several versions of this property. 

\begin{definition} 
Let $G$ be a topological group. \\ 
(1) The group $G$  is called bounded if for any open $V$ containing $1$ there is 
a finite set $F\subseteq G$ and a natural number $k>0$  such that $G=FV^k$. \\ 
(2) The group $G$ is Roelcke bounded if for any open $V$ containing $1$ there is 
a finite set $F\subseteq G$ and a natural number $k>0$  such that $G=V^k FV^k$. \\ 
(3) The group $G$ is Roelcke precompact if for any open $V$ containing 
$1$ there is a finite set $F\subseteq G$ such that $G=VFV$. \\ 
(4) The group $G$ has property ${\bf (OB)_k}$ if for any open symmetric 
$V\not=\emptyset$ there is a finite set $F\subseteq G$ such that $G=(FV)^k$. 
\end{definition} 
It is known that for Polish groups property ${\bf OB}$ is equivalent 
to the property that for any open symmetric $V\not=\emptyset$ there is 
a finite set $F\subseteq G$ and a natural number $k$ such that $G=(FV)^k$.  
Thus when $G$ is non-{\bf OB}, there is an non-empty open $V$ such that 
for any finite $F$ and a natural number $k$, $G\not=(FV)^k$. 
Note that for such $F$ and $k$ there is a real number $\varepsilon$ 
such that some $g\in G$ is $\varepsilon$-distant from $(FV)^k$. 
Indeed, otherwise $(FV)^k V$ would cover $G$.  
 
This explains why in order to  define a suitable class which is 
complementary to {\bf OB} we consider the following property. 

\begin{definition} 
A metric group $G$ is called uniformly non-{\bf OB} if there is an open 
symmetric $V\not=\emptyset$ so that for any natural numbers $m$ and $k$ 
there is a real number $\varepsilon$ such that for any $m$-element subset 
$F\subset G$ there is $g\in G$ which is $\varepsilon$-distant from $(FV)^k$. 

Uniform non-boundedness, uniform non-Roelcke boundedness, uniform non-Roelcke 
precompactness and uniform non-${\bf (OB)_k}$ are defined by the same scheme. 
\end{definition}  

It is clear that in the case discrete groups if a symmetric 
subset $V$ has the property that $G\not= (FV)^k$ for all finite 
$F\subset G$ and natural numbers $k$, then the corresponding 
uniform version also holds.

\begin{proposition} \label{uniform} 
The following classes of metric groups are bountiful: \\ 
(1) The class of uniformly non-bounded groups; \\ 
(2) The class of uniformly non-Roelcke bounded groups; \\  
(3) The class of uniformly non-Roelcke precompact groups; \\ 
(4) The class of uniformly non-${\bf (OB)_k}$-groups; \\   
(5) The class of uniformly non-${\bf (OB)}$-groups.  
\end{proposition} 
 
{\em Proof.} 
Let us consider the class of uniformly non-${\bf (OB)}$-groups. 
Let $\mathcal{K}_{0}$ be the class of all continuous 
metric structures $\langle G, P,Q \rangle$ with the axioms 
stating that $G$ is a group and $P:G\rightarrow [0,1]$ and   
$Q:G\rightarrow [0,1]$ are  unary predicates on $G$ with $Q(1)=0$ so that  
$$ 
sup_x min (P(x),Q(x))= sup_x |P(x)-P(x^{-1})|=0 \mbox{ and } inf_x |P(x)-1/2| =0,  
$$
$$ 
sup_x |Q(x)-Q(x^{-1})|=0 \mbox{ and } inf_x |Q(x)-1/2| =0,  
$$
$$
\mbox{ and for all rational }\varepsilon \in [0,1]  
$$ 
$$
sup_x  min ( \varepsilon \dot{-} Q(x), inf_y (max(d(x,y)\dot{-} 2\varepsilon , \varepsilon \dot{-}  P(y)))=0.   
$$ 
Note that the last axiom implies that any neighbourhood of an element from the nullset of $Q$ 
contains an element with non-zero $P$. 

For any natural $m$ and $k$ and any rational $\varepsilon$ let us consider the following condition (say $\theta (m,k,\varepsilon )$): 
$$ 
sup_{x_1 ...x_m }  inf _{x}   sup_{y_1 ...y_k} min(P(y_1),...,P(y_n), (\varepsilon \dot{-} min_{w\in W_{m,k}}(d(x,  w))))=0, 
$$ 
$$ 
\mbox{ where } W_{m,k} \mbox{ consists of all words of the form } x_{i_1} y_1 x_{i_2} y_2 ...x_{i_k}y_k . 
$$ 
If $G$ is a uniformly non-${\bf (OB)}$-group, then find an open symmetric $V$ 
such that for any natural numbers $m$ and $k$ 
there is a real number $\varepsilon$ such that for any $m$-element subset 
$F\subset G$ there is $g\in G$ which is $\varepsilon$-distant from $(FV)^k$. 
We interpret $Q(x)$ by $d(x,V)$ and $P(x)$ by $d(x, G\setminus V)$ 
(possibly normalizing them to satisfy the axioms of $\mathcal{K}_0$). 
Then observe that $\langle G,P,Q\rangle\in \mathcal{K}_0$ and for any 
natural numbers $m$ and $k$ there is a rational number $\varepsilon$ 
so that $\theta(m,k,\varepsilon )$ holds in $(G,P,Q)$.  
By the L\"{o}wenheim-Skolem theorem for continuous logic  
any infinite subset $C$ of such a structure is contained in an elementary 
submodel of the same density character as $C$. 
To verify uniform {\bf non-(OB)} in such a submodel take the complement 
of the nullset of $P(x)$ as an open symmetric subset. 
This proves statement (5).  

All remaining cases are considered in a similar way. 
$\Box$ 
\bigskip


\subsection{Unbounded actions} 
\paragraph{{\bf Negation of (T).}}

Let a topological group $G$ have a strongly continuous unitary 
representation on a Hilbert space ${\bf H}$.  
A closed subset $Q\subset G$ 
has an {\em almost $\varepsilon$-invariant unit vector} in ${\bf H}$ if 
$$ 
\mbox{ there exists }v\in {\bf H} \mbox{ such that } 
sup_{x\in Q} \parallel x\circ v - v\parallel < \varepsilon
\mbox{ and } \parallel v\parallel =1.  
$$ 
We call a closed subset $Q$ of the group $G$ a {\em Kazhdan set} 
if there is $\varepsilon$ with the following property: 
every unitary representation of $G$ on a Hilbert space 
with almost $(Q,\varepsilon )$-invariant unit vectors also has 
a non-zero invariant vector.  
If the group $G$ has a compact Kahdan subset then it is said that 
$G$ has property ${\bf (T)}$ of Kazhdan. 

If we want to consider unitary representations in 
continuous logic we should fix continuity moduli 
for the corresponding binary functions 
$G\times B_n \rightarrow B_n$ induced by the action, 
where $B_n$ is the $n$-ball of the corresponding Hilbert space. 
In fact if $G$ is $\sigma$-locally compact, 
then we can present $G$ as the union of a chain 
of compact subsets $K_1 \subseteq K_2 \subseteq ...$ 
and consider continuity moduli for the corresponding 
functions $K_m \times B_n \rightarrow B_n$. 
Note that each $B_k$ and $K_l$ will be considered as sorts 
of a continuous structure. 
In this version of continuous logic we do not assume that 
the diameter of a sort is bounded by 1. 
It can become any rational number.  

We can now slightly modify the definition of a Kazhdan set 
as follows. 

\begin{definition} 
Let $G$ be the union of a chain of closed subsets 
$K_1 \subseteq K_2 \subseteq ...$  of bounded diameters.  
Let $\mathcal{F} = \{ F_1 , F_2 , ... \}$ be a family of continuity 
moduli for continuous function $K_i \times B_i \rightarrow B_i$. 

We call a closed subset $Q$ of the group $G$ 
an $\mathcal{F} $-Kazhdan set if there is $\varepsilon$ 
with the following property: 
every $\mathcal{F} $-continuous unitary representation of $G$ 
on a Hilbert space with almost $(Q,\varepsilon )$-invariant unit vectors 
also has a non-zero invariant vector.  
\end{definition} 

Let us consider such actions in continuous logic. 
We treat a Hilbert space over $\mathbb{R}$ 
exactly as in Section 15 of \cite{BYBHU}. 
We identify it with a many-sorted metric structure 
$$
(\{ B_n\}_{n\in \omega} ,0,\{ I_{mn} \}_{m<n} ,
\{ \lambda_r \}_{r\in\mathbb{R}}, +,-,\langle \rangle ),
$$
where $B_n$ is the ball of elements of norm $\le n$, 
$I_{mn}: B_m\rightarrow B_n$ is the inclusion map, 
$\lambda_{r}: B_m\rightarrow B_{km}$ is scalar 
multiplication by $r$, with $k$ the unique integer 
satisfying $k\ge 1$ and $k-1 \le |r|<k$; 
futhermore, $+,- : B_n \times B_n \rightarrow B_{2n}$ 
are vector addition and subtraction and 
$\langle \rangle : B_n \rightarrow [-n^2 ,n^2 ]$ 
is the predicate of the inner product. 
The metric on each sort is given by 
$d(x,y) =\sqrt{ \langle x-y, x-y \rangle }$.   
For every operation the continuity modulus is standard.  
For example in the case of $\lambda_r$ this is $\frac{z}{|r|}$. 

Stating existence of infinite approximations of orthonormal bases 
(by a countable family of axioms, see Section 15 of \cite{BYBHU}) 
we assume that our Hilbert spaces are infinite dimensional. 
By \cite{BYBHU} they form the class of models of a complete 
theory which is $\kappa$-categorical for all infinite $\kappa$, 
and admits elimination of quantifiers. 

This approach can be naturally extended to complex Hilbert spaces, 
$$
(\{ B_n\}_{n\in \omega} ,0,\{ I_{mn} \}_{m<n} ,
\{ \lambda_c \}_{c\in\mathbb{C}}, +,-,\langle \rangle_{Re} , \langle \rangle_{Im} ). 
$$
We only extend the family 
$\lambda_{r}: B_m\rightarrow B_{km}$, $r\in \mathbb{R}$, 
to a family $\lambda_{c}: B_m\rightarrow B_{km}$, $c\in \mathbb{C}$, 
of scalar products by $c\in\mathbb{C}$, with $k$ 
the unique integer satisfying $k\ge 1$ and $k-1 \le |c|<k$. 

We also introduce $Re$- and $Im$-parts of the inner product. 

If we remove from the signature of complex Hilbert spaces 
all scalar products by $c\in \mathbb{C}\setminus \mathbb{Q}[i]$, 
we obtain a countable subsignature 
$$
(\{ B_n\}_{n\in \omega} ,0,\{ I_{mn} \}_{m<n} ,
\{ \lambda_c \}_{c\in\mathbb{Q}[i]}, +,-,\langle \rangle_{Re} , \langle \rangle_{Im} ),
$$
which is {\em dense} in the original one: \\ 
if we present $c\in \mathbb{C}$ by a sequence $\{ q_i \}$ 
from $\mathbb{Q}[i]$ converging to $c$, 
then the choice of the continuity moduli of 
the restricted signature still guarantees that 
in any sort $B_n$ the functions $\lambda_{q_i}$ 
form a sequence which converges to  $\lambda_c$ 
with respect to the metric 
$$ 
sup_{x\in B_n} \{ |f^M (x) - g^M (x)| : M \mbox{ is an $L$-structure } \}.  
$$  
This obviously implies that the original language 
of Hilbert spaces is {\em separable}. 
In particular we may apply Remark \ref{FaHaSh}.

\bigskip

Let us consider a class of continuous metric structures 
which are unions of the many-sorted structures 
$$ 
(\{ K_n \}_{n\in\omega}, \cdot , ^{-1}, 1 ),  
$$
corresponding to groups $G$ presented as $\bigcup_{n\in\omega} K_n$, 
together with metric structures of complex Hilbert spaces  
$$ 
(\{ B_n\}_{n\in \omega} ,0,\{ I_{mn} \}_{m<n} ,
\{ \lambda_c \}_{c\in\mathbb{C}}, +,-,\langle \rangle_{Re} , \langle \rangle_{Im} ).
$$
The operation $\cdot$ (and $^{-1}$) is considered as 
a family of maps $K_n \times K_n \rightarrow K_{n+1}$. 
(maps $K_n\rightarrow K_{n+1}$ respectively), $n\in\omega$.  
Such a structure (say $A(G,{\bf H})$) also contains 
a binary operation $\circ$ of an action which is defined 
by a family of appropriate maps  
$K_n \times B_m \rightarrow B_{m}$. 
When we add the obvious continuous $sup$-axioms that 
the action is linear and unitary, we obtain an axiomatizable 
class $\mathcal{K}_{GH}$. 
We do not state exactly which continuity moduli would 
correspond to these operations. 
In fact this depends on groups and actions we want 
to have in  $\mathcal{K}_{GH}$.  
By $\mathcal{K}_{GH}(\mathcal{F} )$ we denote the 
corresponding class with continuity moduli $\mathcal{F} $. 

Assuming that continuity moduli $\mathcal{F} $ are fixed 
let $\mathcal{K}_{aiv}(\mathcal{F} )$ be the subclass of 
$\mathcal{K}_{GH}$ axiomatizable by the axioms 
$$ 
inf_{v\in B_{m}} sup_{x\in K_{n}} 
max(\parallel x\circ v - v\parallel \dot- \mbox{ } {1\over n}, |1- \parallel v\parallel \mbox{ } |)=0
\mbox{ , } m,n\in\omega \setminus \{ 0\}, 
$$ 
which in fact say that each $K_n$ has an almost 
${1\over n}$-invariant unit vector in ${\bf H}$.

Below we will only consider metric groups which have presentations 
$G=\bigcup_{i\in\omega} K_i$, where $\{ K_i :i\in\omega\}$ 
is an increasing sequence of closed  subsets of diameters 
$d_1 \le d_2 \le ... $ with 
$\{ 1 \} \cup K_n \cdot K_n \cup K^{-1}_n \subseteq K_{n+1}$. 
Moreover we will assume that every function 
$K_n\times K_n \rightarrow K_{n+1}$ induced by the 
multiplication is uniformly continuous with respect to 
some fixed family of continuity moduli $\mathcal{F}_0$.  
When $G$ is a $\sigma$-locally compact group 
then such $\mathcal{F}_0$ and such a decomposition obviously exist.

\begin{proposition} \label{nT} 
Let $G$ be a metric group having a presentation 
$G= \bigcup_{i\in\omega} K_i$ into an icreasing 
sequence of closed subsets of bounded diameters as above, 
so that no  $K_i$ is an $\mathcal{F} $-Kazhdan set for $G$. 
Then there is a metric structure in $\mathcal{K}_{aiv}(\mathcal{F} )$  
which naturally expands 
$(\{ K_n \}_{n\in\omega}, \cdot , ^{-1}, 1 )$ 
so that $G$ does not have non-zero fixed vectors. 
In particular any $\sigma$-locally compact group $G$ 
without property {\bf (T)} for $\mathcal{F} $-actions 
has such an expansion to ${\bf H}$. 
\end{proposition} 

{\em Proof.} 
Consider the particular case of the proposition. 
Let $G=\bigcup_{i\in\omega} K_i$, where $\{ K_i :i\in\omega\}$ 
is an increasing sequence of compact neighborhoods of $1$ 
with $K_n \cdot K_n \cup K^{-1}_n \subseteq K_{n+1}$. 
We know that for any natural $n$ and rational $0<q<1$ there is 
a unitary $\mathcal{F}$-continuous representation of $G$ 
on a Hilbert space ${\bf H}$ which has an almost $(K_n ,q)$-invariant 
unit vector but does not have a non-zero invariant vector. 
Decomposing a basis of ${\bf H}$ into a union of an infinite family 
of pairwise disjoint infinite subsets  (labelled by pairs $(n,q)$) 
and defining representations above on the corresponding subspaces, 
we find an unitary $\mathcal{F}$-representation of $G$ on ${\bf H}$ 
without non-zero invariant vectors so that for any natural $n$ and rational 
$q<1$ there is an almost $(K_n ,q)$-invariant unit vector in ${\bf H}$. 

Let us define a required metric structure 
(denoted by $A(G,{\bf H})$ as above) 
as a union of the many-sorted structure 
$$ 
(\{ K_n \}_{n\in\omega}, \cdot , ^{-1}, 1 ),  
$$
corresponding to the group $G$, together with 
the metric structure of the separable Hilbert space 
$$ 
(\{ B_n\}_{n\in \omega} ,0,\{ I_{mn} \}_{m<n} ,
\{ \lambda_c \}_{c\in\mathbb{C}}, +,-,\langle \rangle ).
$$
The operation $\cdot$ (and $^{-1}$) is considered as 
a family of maps $K_n \times K_n \rightarrow K_{n+1}$. 
(maps $K_n\rightarrow K_{n+1}$ respectively), $n\in\omega$.  
The structure $A(G,{\bf H})$ also contains a binary 
operation $\circ$ of an action which is defined by 
a family of appropriate maps  
$K_n \times B_m \rightarrow B_m$. 
Since the representation is unitary, any element 
of $G$ preserves each $B_m$.  
The axioms of $\mathcal{K}_{inv}(\mathcal{F})$ obviously 
hold in $A(G,{\bf H})$. 

The argument in the general case is the same. 
$\Box$ 

\bigskip 

\begin{remark} 
{\em When any compact subset of $G$ is 
contained in some $K_i$, the group $G$ from 
the formulation above does not have property {\bf (T)} 
for $\mathcal{F}$-continuous representations. 
For example this happens when each $K_i$ 
is the closed $d_i$-ball of $1$. }
\end{remark}
\bigskip 

The following corollary can be considered 
as a kind of bountifulness for groups without {\bf (T)}. 
It follows from the proposition and remark above and 
the L\"{o}wenheim-Skolem theorem for continuous logic 
(Remark \ref{FaHaSh}). 

\begin{corollary} 
Let $G$ be a metric group having a presentation 
$G= \bigcup_{i\in\omega} K_i$ into an icreasing 
sequence of closed balls of $1$ of bounded diameters, 
so that no  $K_i$ is an $\mathcal{F}$-Kazhdan set for $G$. 
Then for any infinite subset $C\subset G$ there is a closed 
elementary (in continuous logic) subgroup of $G$ containing $C$, 
with the same density character as $C$ 
and without property {\bf (T)} for 
$\mathcal{F}$-continuous representations. 
\end{corollary} 


\paragraph{{\bf Non-FH-actions.}} 

To consider non-{\bf FH} let us fix a binary function 
$\nu :\omega\times\omega \rightarrow\omega$ 
which is increasing in each argument. 
We now define $\mathcal{K}_{GH}(\nu )$, a class 
of continuous metric structures which are unions 
of many-sorted $\mathcal{F}_0$-continuous structures 
$$ 
(\{ K_n \}_{n\in\omega}, \cdot , ^{-1}, 1 ),  
$$
corresponding to groups $G$ presented as 
$\bigcup_{n\in\omega}K_n$ (with assumptions as before Proposition \ref{nT}), 
together with metric structures of real Hilbert spaces  
$$ 
(\{ B_n\}_{n\in \omega} ,0,\{ I_{mn} \}_{m<n} ,
\{ \lambda_c \}_{r\in\mathbb{R}}, +,-,\langle \rangle ).
$$   
We also add a binary $\mathcal{F}$-continuous operation 
$\circ$ of an action defined by a family of appropriate 
maps $K_n \times B_m \rightarrow B_{\nu (n,m)}$. 
Obvious continuous $sup$-axioms that the action 
is isometric give an axiomatizable 
class $\mathcal{K}_{GH}(\nu ,\mathcal{F} )$.

The following statement is a straightforward application of the continuous L\"{o}wenheim-Skolem theorem. 

\begin{proposition} 
Let $G$ be a metric group having a presentation  
$G= \bigcup_{i\in\omega} K_i$ into an icreasing 
sequence of closed subsets of bounded diameters, 
so that there is an isometric action of $G$ on 
a real Hilbert space inducing a structure from 
$\mathcal{K}_{GH}(\nu, \mathcal{F} )$ without fixed points. 
Then for any infinite subset $C\subset G$ 
there is a closed continuously elementary subgroup of $G$ 
containing $C$, with the same density character as $C$ 
and without property {\bf FH} for $\mathcal{F}$-actions. 
\end{proposition} 

Note that any $\sigma$-locally compact metric group 
$G$ without property {\bf FH} for $\mathcal{F}$-actions 
belongs to the class $K_{GH} (\nu ,\mathcal{F})$ 
(for appropriate $\nu$ and $\mathcal{F}_0$). 
We will now see this in a slightly stronger form. 
Fix a function $\eta$ which assigns to a natural number $k$ 
a pair $(l,s)\in \omega\times\omega$.  
To obtain $\mathcal{K}_{nFH}(\nu ,\eta ,\mathcal{F})$ take the 
subclass of $\mathcal{K}_{GH}(\nu )$ axiomatizable by the axioms 
$$ 
sup_{v\in B_{k}} inf_{x\in K_{l}} 
({1\over s} \dot{-}  \parallel x\circ v - v\parallel  )=0  
\mbox{ , for } \eta (k) =(l,s)   
$$ 
(saying that each vector of $B_k$ is moved by some  
element of $K_l$ by approximately $ {1\over s}$). 

\begin{proposition} 
For any $\sigma$-locally compact metric group $G$ 
without property {\bf FH} for $\mathcal{F}$-continuous 
actions there are functions $\nu$ and $\eta$ so that 
$G$ belongs to the class $\mathcal{K}_{nFH}(\nu ,\eta ,\mathcal{F})$. 
\end{proposition} 

{\em Proof.} 
Let $G=\bigcup_{i\in\omega} K_i$, where $\{ K_i :i\in\omega\}$ 
is an increasing sequence of compact neighborhoods of $1$ 
with $K_n \cdot K_n \cup K^{-1}_n \subseteq K_{n+1}$. 
Fix appropriate $\mathcal{F}_0$. 
We know that there is an affine $\mathcal{F}$-continuous 
isometric action of $G$ on a real Hilbert space ${\bf H}$ 
which does not have a fixed vector. 
Since for every $v\in {\bf H}$ the map 
$$
G \rightarrow {\bf H} \mbox{ , } g \rightarrow gv 
$$ 
is continuous there are $k\in\omega$ and an open subset 
of $G$ which maps $0$ into $B_k$. 
In particular for any $n$ we can find such a $k$ so that 
all elements of $K_n$ map $0$ into $B_k$. 
This means that $K_n \circ B_m \subseteq B_{k+m}$. 
This defines a function 
$\nu :\omega\times\omega \rightarrow\omega$ so that 
the continuous structure corresponding to the action 
belongs to $\mathcal{K}_{GH}(\nu ,\mathcal{F} )$. 

Let us show that this structure belongs to 
$\mathcal{K}_{nFH}(\nu ,\eta ,\mathcal{F})$ 
for appropriate $\eta$. 
Since $G$ does not fix any point,  each orbit of $G$ 
is unbounded (Proposition 2.2.9 of \cite{BHV}). 
Thus there is $g\in G$ so that $B_k \cap g (B_k) =\emptyset$. 
In particular there is $s\in \mathbb{N}$ such that 
$ {1\over s} \le \parallel g\circ v - v\parallel$ 
for all $v\in B_k$.  
We define $\eta (k)$ to be $(l,s)$, where $l$ is chosen so 
that $K_l$ contains $g$ as above. 
$\Box$ 

\bigskip 

\paragraph{{\bf Non-F$\mathbb{R}$-actions.}}   

To consider non-{\bf F}$\mathbb{R}$ we apply similar ideas. 
Let us fix a binary function 
$\nu :\omega\times\omega \rightarrow\omega$ 
which is increasing in each argument.  
We now define $\mathcal{K}_{GR}(\nu ,\mathcal{F} )$, a class 
of continuous metric structures 
which are unions of the many-sorted structures 
$$ 
(\{ K_n \}_{n\in\omega}, \cdot , ^{-1}, 1 ),  
$$
corresponding to groups $G$ presented as $\bigcup_{n\in\omega} K_n$, 
together with metric structures of pointed real trees  
$$ 
(\{ B_n\}_{n\in \omega} ,0,\{ I_{mn} \}_{m<n} ), 
$$   
where $B_n$ is the $n$-ball of $0$ in the tree. 
It is shown in \cite{carlisle} 
that the class of pointed real trees is axiomatizable 
in contimuous logic by axioms of $0$-hyperbolicity and 
the approximate midpoint property. 
 
We also add a binary operation $\circ$ of an action 
defined by a family of appropriate $\mathcal{F}$-continuous 
maps $K_n \times B_m \rightarrow B_{\nu (n,m)}$. 
Obvious continuous $sup$-axioms that the action 
is isometric give an axiomatizable 
class $\mathcal{K}_{GR}(\nu ,\mathcal{F})$. 

As in the non-{\bf FH}-case we have the following straightforward statement. 

\begin{proposition} \label{211} 
Let $G$ be an $\mathcal{F}_0$-continuous group with respect 
to a presentation $G= \bigcup_{i\in\omega} K_i$ into 
an icreasing sequence of closed subsets of bounded diameters, 
so that there is an isometric action of $G$ on a real tree 
inducing a structure from $\mathcal{K}_{GR}(\nu ,\mathcal{F})$ 
without fixed points. 
Then for any infinite subset $C\subset G$ there 
is a closed continuously elementary subgroup of $G$ 
containing $C$, with the same density character as $C$ 
and without property {\bf F}$\mathbb{R}$ for 
$\mathcal{F}$-continuous actions. 
\end{proposition} 

Fix a function $\eta$ which assigns to a natural 
number $k$ a pair $(l,s)\in \omega\times\omega$.  
To obtain $\mathcal{K}_{nFR}(\nu ,\eta ,\mathcal{F})$ take 
the subclass of $\mathcal{K}_{GR}(\nu )$ axiomatised  
by the axioms 
$$ 
sup_{v\in B_{k}} inf_{x\in K_{l}} 
(  {1\over s} \dot{-} d(x\circ v,v  ))=0  \mbox{ , for } \eta (k)=(l,s)   
$$ 
(saying that each element of $B_k$ is moved by some  
element of $K_l$ by approximately $ {1\over s}$).

\begin{proposition} 
For any $\sigma$-locally compact metric 
group $G$ without property {\bf F}$\mathbb{R}$ 
for strongly $\mathcal{F}$-continuous actions there are functions 
$\nu$ and $\eta$ so that $G$ belongs to the class 
$\mathcal{K}_{nFR}(\nu ,\eta ,\mathcal{F})$. 
\end{proposition} 

{\em Proof.} 
Let $G=\bigcup_{i\in\omega} K_i$, where $\{ K_i :i\in\omega\}$ 
is an increasing sequence of compact neighborhoods of $1$ 
with $K_n \cdot K_n \cup K^{-1}_n \subseteq K_{n+1}$. 
Find appropriate moduli $\mathcal{F}_0$ making $G$ 
an $\mathcal{F}_0$-continuous structure. 
We know that there is an $\mathcal{F}$-continuous 
isometric action of $G$ on a real tree $T$ which 
does not have a fixed point. 
Since for every $v\in T$ the map 
$$
G \rightarrow T \mbox{ , } g \rightarrow gx 
$$ 
is continuous there are $k\in\omega$ and an open subset 
of $G$ which maps $0$ into $B_k$. 
In particular for any $n$ we can find a $k$ so that 
all elements of $K_n$ map $0$ into $B_k$. 
This means that $K_n \circ B_m \subseteq B_{k+m}$. 
This defines a function 
$\nu :\omega\times\omega \rightarrow\omega$ so that 
the continuous structure corresponding to the action 
belongs to $\mathcal{K}_{GR}(\nu ,\mathcal{F} )$. 

Let us show that this structure belongs to 
$\mathcal{K}_{nFR}(\nu ,\eta ,\mathcal{F})$ for appropriate $\eta$. 
If $G$ has a hyperbolic element $g$ of hyperbolic length $r$ 
(i.e. there is a line $L$ so that $g$ $r$-shifts all points of $L$), 
then $\eta$ is constant where $l$  is chosen so 
that $K_l$ contains this hyperbolic element 
and $s$ is chosen so that $ {1\over s} \le r$. 

Consider the case when $G$ consists of elliptic elements (i.e. fixing points). 
Since $G$ does not fix any point, by a well-known 
argument $G$ fixes an end (\cite{serre}, Section 6.5, Exercise 2). 
Let $L_0$ be the half-line starting from $0$ which 
represents this end and let $v_1 ,....,v_i ,....$  
be a cofinal $\omega$-sequence in $L_0$ with $d(v_i ,v_{i+1})\ge 1$. 
Then we may assume that $G$ is the union of a strictly 
increasing chain of stabilizers $G_i$ of $v_i$. 
Since the action is continuous, all $G_i$ are closed. 

Having $k$ find $j$ with $v_{j-1} \not\in B_k$ (thus $v_j \not\in B_k$). 
Since any arc linking $v_j$ with an element from $B_k$ must 
contain $v_{j-1}$, we see that if $g\in G_j$ fixes a point 
of $B_k$ then it fixes $v_{j-1}$. 
In particular $G_j$ does not fix any element of $B_k$. 
Since $d(v_j ,v_{j-1}) \ge 1$ any point of $B_k$ can 
be taken by some element of $G_j$ at a distance greater than 1.   
Thus to define $\eta (k)=(l,s)$, we choose $l$ so that $K_l$ contains 
an element of $G_j$ not fixing $v_{j-1}$. 
We define $s=1$. 
$\Box$

\bigskip

\section{Separably categorical locally compact groups} 

If $G$ is a locally compact group, then $G$ admits 
a compatible complete left invariant metric $d(x,y)$ 
(\cite{becker}, 3.C.2). 
We may assume that $d(x,y)$ satisfies $d(x,y) < 1$ 
(it can be replaced by $\frac{d(x,y)}{d(x,y)+1}$). 
We thus may consider locally compact groups in continuous logic 
as a class of continuous metric structures  
$$ 
(G, d, \cdot , ^{-1}, 1 ),  
$$ 
together with fixed continuity moduli for functional symbols.

\begin{theorem} \label{catlc}
Let $G$ be a separably categorical locally compact 
non-compact group. 
Then there is a compact clopen subgroup $H<G$ 
which is invariant with respect to all metric automorphisms of $G$, 
and the induced action of $Aut(G,d)$ on the coset space $G/H$ 
is oligomorphic. 

If the connected component of the unity $G^0$ is not trivial, 
$H$ can be taken to be $G^0$. 
In this case and in the case when $d$ is two-sided-invariant, 
the subgroup $H$ is normal and $G/H$ is an $\omega$-categorical 
discrete group. 
\end{theorem} 

We start with the following preliminaries. 
We may assume that $G$ is not discrete. 
There is a rational number $\rho <1$ such that the $\rho$-ball 
of the unity $B_{\rho} (1)=\{ x\in G: d(x,1)\le\rho\}$ is compact. 
In particular $B_{\rho}(1)$ is a subset of $acl(\emptyset )$ in $G$ 
(the condition $d(x,1)\le \rho$ defines a totally bounded, complete subset in any elementary extension of $G$). 
Thus any $B^n_{\rho}(1)$ also is a subset of $acl(\emptyset )$. 
Let $G_{\rho}$ be the subgroup generated by $B_{\rho}(1)$. 
Note that for any $g\in G_{\rho}$ the open ball 
$$
B_{<\rho}( g)=\{ x\in G: d(x,g)<\rho\}= \{ x\in G: d(g^{-1} x,1) <\rho \}
$$
is a subset of $G_{\rho}$; thus $G_{\rho}$ is an open (in fact clopen) subgroup. 
If $G$ has a non-trivial connected component of the unity $G^0$, 
we may assume that $G_{\rho}=G^0$. 
Note that when $d$ is a two-sided-invariant metric, $G_{\rho}$ 
is a normal subgroup of $G$. 
 
\begin{lemma} 
Assume that $G$ is a separably categorical 
locally compact group. 
Then under the circumstances above 
the predicate $P(x)=d(x,G_{\rho})$ is definable in $G$. 
\end{lemma} 

{\em Proof.} 
Since the space $(S_n(T),d)$ is compact, 
for every $\varepsilon$ there is 
a finite set of types which is $\varepsilon$-dense 
in the set of types of elements of $\bigcup_{n>0} B^{n}_{\rho}(1)$. 
Thus there is a number $n$ such that 
the $\varepsilon$-neighbourhood of $B^n_{\rho}(1)$ 
contains the zeroset of  $P(x) = dist (x,G_{\rho})$.
If $(N,Q)$ is an elementary extension of $(G,P)$ then 
$(N,Q)$ satisfies the condition 
$$ 
sup_x inf_{y_1} ...inf_{y_n}max(d(y_1 ,1)\dot{-}\rho ,...,  d(y_n ,1)\dot{-}\rho, |Q(x) -d(x,y_1 \cdot ...\cdot y_n )|\dot{-}\varepsilon )=0, 
$$ 
i.e. the $\varepsilon$-neighbourhood of $B^n_{\rho}(1)$ 
contains the zeroset of $Q(x)$.   
In particular the zeroset of $Q$ coincides with the closure of $G_{\rho}$, 
i.e is $G_{\rho}$ itself and is a subset of $acl(\emptyset )$. 
Since $(N,Q)$ is an elementary extension of $(G,P)$,
$Q(x)$ is the distance from the zeroset of $Q$ (see Theorem 9.12 in \cite{BYBHU}). 
In particular any automorphism of $N$ preserves $Q$. 
Using Corollary 9.11 of \cite{BYBHU} (cited in  Introduction above)   
we see that  $P(x)$ is a definable predicate. 
$\Box$

\bigskip

\begin{lemma} 
Under the circumstances above  
there is a natural number $n$ so that 
$G_{\rho}=B^n_{\rho}(1)$. 
In particular $G_{\rho}$ is compact. 
\end{lemma} 

{\em Proof.} 
If $G_{\rho}\not=B^n_{\rho}(1)$ for all $n\in\omega$, there are 
positive rational numbers $\varepsilon_1 ,...,\varepsilon_n ,...$ so that 
the $\varepsilon_n$-neighbourhood of $B^n_{\rho}(1)$ does not cover $G_{\rho}$.    
Thus all statements 
$$
sup_{x_1 ...x_n}(min(\varepsilon_n \dot{-} d(x,x_1 \cdot ....\cdot x_n ), \rho\dot{-}d(1,x_1 ),..., \rho\dot{-}d(1,x_n )))=0
$$ 
are finitely consistent together with $P(x)=0$. 
By compactness of continuous logic we obtain a contradiction. 
$\Box$  

\bigskip 

Since $G_{\rho}$ is a characteristic subgroup of $G$ with respect to 
the automorphism group of the metric structure $G$, 
we see that $Aut(G,d)$ acts correctly on $G/G_{\rho}$ 
by permutations of $G/G_{\rho}$. 
Note that $G/G_{\rho}$ is a discrete space with respect to 
the topology induced by the topology of $G$. 

\begin{lemma} 
The action of $Aut(G,d)$ on $G/G_{\rho}$ is oligomorphic. 
\end{lemma} 

{\em Proof.} 
Since $(G,d)$ is separably categorical, $Aut(G,d)$ 
is approximately oligomorphic on $(G,d)$. 
Thus for every $n$ there is a finite set $F$ of 
$n$-tuples from $G$ such that the set of orbits 
meeting $F$ is $\rho$-dense in $(G,d)$. 
In particular for any $g_1 ,...,g_n \in G$ 
there is a tuple $(h_1 ,...,h_n )\in F$ and 
an automorphism $\alpha \in Aut(G,d)$ 
such that $g^{-1}_i \alpha (h_i )\in G_{\rho}$ 
for all $i\le n$. 
$\Box$ 

\bigskip 

To see that Theorem \ref{catlc} follows from lemmas above 
just take $H$ to be $G_{\rho}$.


\bigskip

Institute of Mathematics, \parskip0pt

Wroc{\l}aw University, \parskip0pt

pl.Grunwaldzki 2/4, \parskip0pt

50-384 Wroc{\l}aw, \parskip0pt

POLAND. \parskip0pt

{\em ivanov@math.uni.wroc.pl} 

\end{document}